\documentclass[a4paper]{amsart}
\pdfoutput=1
\usepackage{graphicx}
\usepackage{mathptmx}      
\usepackage{latexsym}
\usepackage{amsmath,amssymb}
\usepackage{color}
\usepackage{enumerate}

\newcommand{\R}{\mathbb R}

\newcommand{\pa}{\mathop{}\!\partial}

\newcommand{\eps}{\varepsilon}
\newcommand{\ve}[1]{{\mathbf {#1}}}

\newtheorem{theorem}{Theorem}[section]

\begin{document}

\title{An invariant region for the collisional dynamics of 
two bodies on Keplerian orbits}

\maketitle

\author{Dario Benedetto}

\address{Dipartimento di Matematica, {\it Sapienza} Universit\`a di Roma }
\email{benedetto@mat.uniroma1.it}

\author{Flavia Lenti}

\address{Dipartimento di Scienza e di alta Tecnologia
Universit\`a dell'Insubria}
\email{flavia.lenti@uninsubria.it}

\begin{abstract}
We study the dynamics of two bodies moving on
elliptic Keplerian orbits around a fixed center 
of attraction and interacting
only by means of elastic or inelastic collisions. We show that
there exists a bounded invariant region: for
suitable values of the total energy and the total angular
momentum (explicitly computable) the orbits of the bodies remain elliptic,
whatever are the number and the details 
of the collisions.  We show that there exists a
bounded invariant region even in the case of two bodies interacting by
short range potential.
\end{abstract}



\section{Introduction}
\label{sezione:introduzione}

The interest in the collisional dynamics in a
planetary system goes back to Poincar\'e. In particular, 
in \cite{poincare} 
he studies
the planetary three-body problem, with 
one body of large mass (the Sun)
and two bodies of small mass (the planets).
He indicates how to find periodic 
solutions (of {\it deuxi\`eme esp\`ece}),  
as perturbations of  periodic collisional solutions
he can construct when the mass of the planets are infinitely small
and the distance between them becomes infinitely small.
In this approximation, 
the two bodies are on Keplerian ellipses until the ``choc'' (i.e.
the interaction), which moves the bodies on two other Keplerian 
ellipses. During the interaction, only the total energy and the total
momentum are conserved, and the choc acts 
as an elastic collision. It should be noted that
the collision can moves the bodies also on hyperbolic orbits but
Poincar\'e is only interested in elliptic case.

In this work, we prove that, for this collisional dynamics,
there exists an invariant bounded region of
of positive measure in the phase space.
More precisely, we consider two
bodies moving on elliptic Keplerian orbits around a center,
interacting only by means of collisions. A collision changes the
orbital parameters of the bodies, and a sequence of collisions can
move one of the bodies out of the system, on parabolic or hyperbolic
orbits.  We show that for suitable values of the total energy and
the total angular momentum (easily computable), 
the bodies remain on elliptic colliding orbits. 
Moreover, we extend this result to the case of two point particles
interacting by means a bounded short range potential.

We are  neglecting the gravitational interaction between the bodies 
and the influence on the attractive center.
These approximations can be only justified for very light
particles and for a few revolution times
(the time between two consecutive
collisions can be very long with respect to the 
revolution period of the particles). 
Considering this,  our result has the following 
interpretation for real systems:
the particles can not leave the system due to the collisions,
unless other perturbations change enough the
orbital parameters. 

Our result seems not known in the literature, despite
its simplicity and despite the great interest in research on 
planetary systems.
We run in the result while we were studying a numerical models
for the dynamics of inelastic particles of planetary rings.
It is  known (see  
\cite{poincare-cosm},  \cite{braich},  and  
\cite{bh})
that the inelasticity of the collisions 
is sufficient to guarantee the persistence of the rings.
In contrast, in the case of  
elastic collisions,  almost all the 
particles leave the system on hyperbolic orbits.  
Here we prove that in the case of only two particles,
the inelasticity is not needed
to avoid that the orbits become hyperbolic or parabolic;
in this sense, two colliding particles are
a stable subsystem of a ring.
This observation and its consequences can be interesting 
for the study of various models of planetary rings.
In particular, in some models it is assumed that the 
collisions are elastic 
for small relative velocities (see \cite{hatzes}
for experimental result on the inelasticity of
ice balls).

From the mathematical point of view, it can be interesting
to study how the particles moves on the invariant region.
Preliminary two-dimensional numerical simulations, 
in which the impact parameter
is randomly chosen, show that 
the two orbits are approximately tangent, for most of the time.
It is difficult to implement a numerical simulation for the more
realistic 
three-dimensional case.
The problem is that it is not easy to find  efficiently
the values of the anomalies which correspond to collisional
configurations.
Some useful suggestion and some technical insight
can be obtained analyzing 
the solutions to the problem 
of finding the critical points of the distance from two 
elliptic orbits (see \cite{gronchi}).

\vskip.3cm
The paper is organized as follows.
In section \ref{sez:problema}, we establish the mathematical notation 
and the exact nature of the problem. 
In section \ref{sezione:invariante}, we analyze a simplified
model considering a two-dimensional {\it dynamics of the orbits}
instead of the dynamics of the to bodies.
We consider more general cases in section \ref{sezione:estensioni} 
(two bodies in $\R^2$) and in section  
 \ref{sezione:dim3}, in which we  
analyze the case of two bodies in $\R^3$; here we also 
analyze the case of point particles interacting by means a
short range potential.

\section{The problem}
\label{sez:problema}

The system we analyze consists in two spherical bodies, of mass $m_1$ and
$m_2$ and radii $R_1$ and $R_2$ respectively,
which are attracted by a fixed gravitational center and interact 
by means of elastic or inelastic collisions.
We indicate with $\ve x_i$ and $\ve v_i$ the position of the center
and the velocity of the body $i=1,2$, 
with $M=m_1+m_2$ the total mass
and with $\mu_i = m_i/M$ the fraction of the total mass
carried by the body $i$.
We choose the units of measure in such that
the gravitational potential energy of the body $i$ 
is $m_i/r_i$,  where $r_i=|\ve x_i|$ is the distance 
from the attracting center.

In the case of hard spheres,  
inelastic collisions can be modeled supposing 
that only a part of the normal impulse is transferred, 
while the tangential one is conserved. Let us denote with 
$\ve n = (\ve x_1 - \ve x_2)/(R_1+R_2)$ the direction 
of the relative position at the moment of the impact, and with $\ve w = \ve v_1 - \ve v_2$ the 
relative velocity. The particles are in the incoming configuration iff $\ve n\cdot \ve w < 0$.
The relative velocity after the collision is
\begin{equation}
\label{w-urto}
\ve w' = (I-\ve n\times \ve n) \ \ve w - ( 1-2\eps) 
(\ve n \times \ve n) \  \ve w
\end{equation}
where $I$ is the identity matrix,  
$(\ve x \times \ve y)_{ij} = x_i y_j$ 
is the tensor product, $\ve n\times \ve n$ is the projector on 
the direction of $\ve n$, $I-\ve n\times \ve n$ is the projector 
on the orthogonal plane to $\ve n$, 
and finally $\eps\in [0,0.5]$ is
the parameter of inelasticity, which is $0$ in the case
of the elastic collision.

The outgoing velocities $\ve v_1'$, $\ve v_2'$
can be obtained from $\ve w'$ using the conservation
of the velocity of the center of mass
$$\ve v' = \mu_1 \ve v_1' + \mu_2 \ve v_2' = \mu_1 \ve v_1 + \mu_2 \ve v_2
= \ve v$$
If the collision is inelastic,
the normal component of $\ve w$ is reduced in modulus:
\begin{equation}
\label{w-normale}
\ve w' \cdot \ve n = -(1-2\eps) \ve w \cdot \ve n,\ \ \ 
\ |\ve w' \cdot \ve n| \le |\ve w \cdot \ve n|,\end{equation}
where $(1-2\eps)\in[0,1)$ is the {\it coefficient of restitution}
which is $1$ in the elastic case.
The kinetic energy 
$T= m_1 |\ve v_1|^2/2 +  m_2 |\ve v_2|^2/2$
decreases and becomes
$$T' = \frac 12 \left(  m_1 |\ve v_1'|^2  +  m_2 |\ve v_2'|^2
 \right) =  T - 2 M \mu_1 \mu_2 \eps (1-\eps) 
(\ve w \cdot \ve n)^2.
$$

We remark that in some model 
it is assumed that also the tangential component
$(I-\ve n \times \ve n)\ \ve w$ is reduced
(see \cite{bh} and \cite{kawai}).
Moreover the restitution coefficient can depend on 
the relative velocity, as shown in a huge number 
of theoretical and experimental studies
(see e.g. \cite{djerassi}, \cite{hertzsch}, and references therein). 
In all these models the kinetic energy decreases.

\vskip3pt

Let us describe what it can be happen to the orbits after a collision
of the two bodies. Before the collision
both energies are 
negative:
$$\frac {m_i}2 \ve v_i^2 - \frac {m_i}{|\ve x_i|} < 0,
\ \ \ i\in 1,2.$$
These conditions are equivalent to 
\begin{equation}
\label{condw}
|\ve v+ \mu_2 \ve w|^2 < 2 |\ve x_1|, \ \ \ 
|\ve v - \mu_1 \ve w|^2 < 2 |\ve x_2|,
\end{equation}
with $|\ve x_1 - \ve x_2|= R_1 + R_2$.
If
\begin{equation}
\label{eq:condw0}
|\ve w| < \min \left( (\sqrt{2|\ve x_1|} - |\ve v|)/\mu_2,
(\sqrt{2|\ve x_2|} - |\ve v|)/\mu_1\right)
\end{equation}
the inequalities \eqref{condw} are satisfied, then, 
after the collision, the orbits  remain elliptic 
because 
$|\ve w'|\le |\ve w|$, as follows from eq.s \eqref{w-urto}, \eqref{w-normale}.
The hypothesis \eqref{eq:condw0} is not sufficient to avoid that
one of the particles leaves the system, because the 
next collision can take place at other points 
with very different values of $\ve w$, $\ve v$,  $\ve x_1$ and $\ve x_2$.
More in general, 
a sequence of collisions on different points
can end with a particle which leaves the system on an hyperbolic orbit.

In the next section  we show how to control the condition of ellipticity,
regardless the collision history.
Let us first note that, for bodies with vanishing radii, 
the possible points of collision
are at most two, which is the maximum number of intersection
of two non identical co-focal Keplerian orbits.
As noted by Poincar\'e in \cite{poincare}, 
if two orbits have two points of intersections
the points and the center of gravity are on a line
or the two orbits are in a plane. 
Then we start our analysis in section \ref{sezione:invariante}
with the two dimensional case, 
also simplifying the dynamics considering two point
particles.

\vskip3pt

\section{The 
invariant region for two point-particles in $\R^2$}
\label{sezione:invariante}

In this section we consider
a simplified bi-dimensional mathematical model. 
We suppose that the bodies are point particles
moving on co-focal Keplerian orbits in a plane.
In order to allow the 
particles collide,  we have to assume 
that the orbits intersect (note that this can happen at most in two points).
Although this condition is satisfied,
the particles can not collide because the cross section
is zero for dimensionless bodies.
In order to avoids this problem, we considered the
{\it dynamics of the orbits}:
we choose as variables the parameters of two orbits $o_1$ and $o_2$,
and we evolve the system with the following procedure:
we choose one of the points of intersection, we consider
in that point
two fictitious particles on the two orbits, 
we choose an impact parameter $\ve n$ and we consider the resulting orbits
 $o_1'$ and $o_2'$
after the collision of the two particles. We will show, in 
Theorem \ref{teo:inter}, that 
for sufficiently low value of
the total energy, 
the new orbits $o_1'$ and $o_2'$ remain intersecting ellipses, 
for any choice between the two 
possible collision points and for any choice of 
the impact parameter $\ve n$. The keys points of the proof are
the conservation of the total momentum,
the decrease of the energy, and the fact that 
the condition of intersection of the orbits is a feature
preserved by the dynamics.

\vskip 3pt
For the orbit $o_i$ of the particle $i$, 
$\omega_i$ is  the angle between the $x$ axis and the 
position of the periapsis (the point of the orbit 
of the minimal distances from the attractive center);
$\vartheta_i$ is the true anomaly, 
 i.e. the angle in the orbital plane between the particle 
and the periapsis 
of the orbit;
$E_i= \ve v_i^2/2 - 1/r_i$ is the specific energy (i.e. the energy 
for unit of mass);
$L_i=r_i^2 \dot \vartheta_i$  is the specific angular momentum;
$e_i = \sqrt{1 + 2 E_i L_i^2}$ is the eccentricity.
The position of the particle $i$ in orbits is given by
\begin{equation}
\label{ri}
\ve x_i = \frac {L_i^2}{1+e_i\cos \vartheta_i} 
\binom {\cos(\vartheta_i + \omega_i)}{\sin(\vartheta_i + \omega_i)}
\end{equation}

\vskip3pt
All the quantities
$E_i$ and $L_i$ are conserved during the Keplerian motion 
but only the combinations $m_1E_1+m_2E_2$
(the total energy) 
and  $m_1L_1+m_2L_2$ (the total angular momentum)
are conserved in the elastic collisions, moreover 
$m_1E_1+m_2E_2$ decreases if the collisions are inelastic.
We fix the initial values of the specific energy
and the specific angular momentum of the whole system:
\begin{equation}
\label{medie}
\begin{array}{l}
L = \mu_1 L_1 + \mu_2 L_2 \\
E = \mu_1 E_1 + \mu_2 E_2 
\end{array}
\end{equation}
We can assume $L\ge 0$ without loss of generality, 
and we consider $E<0$, which is the case  
of a couple of elliptic orbits. 
We rewrite the orbital parameters in terms of the 
differences of the energy and angular momentum:
\begin{equation}
\label{differenze}
\begin{array}{l}
\delta E = E_1 - E_2\\
\delta L = L_1 - L_2
\end{array}
\ \ \ \text{ from which } \ \ \ 
\begin{array}{l}
E_1 = E +\mu_2 \delta E\\
E_2 = E -\mu_1 \delta E\\
\end{array}
\ \ \ \text{ and } \ \ \ 
\begin{array}{l}
L_1 = L +\mu_2 \delta L\\
L_2 = L -\mu_1 \delta L.\\
\end{array}
\end{equation}
From these quantities we can obtain the 
shapes (i.e. the
eccentricities) and the dimensions of the orbits. 
Theirs positions in the framework are specified by the angles $\omega_i$, but 
we are only interested in the relative position of the two orbits,
which is given by $\delta \omega = \omega_1 -\omega_2$.

\vskip3pt



\begin{figure*}[h]
  \includegraphics[scale=0.3]{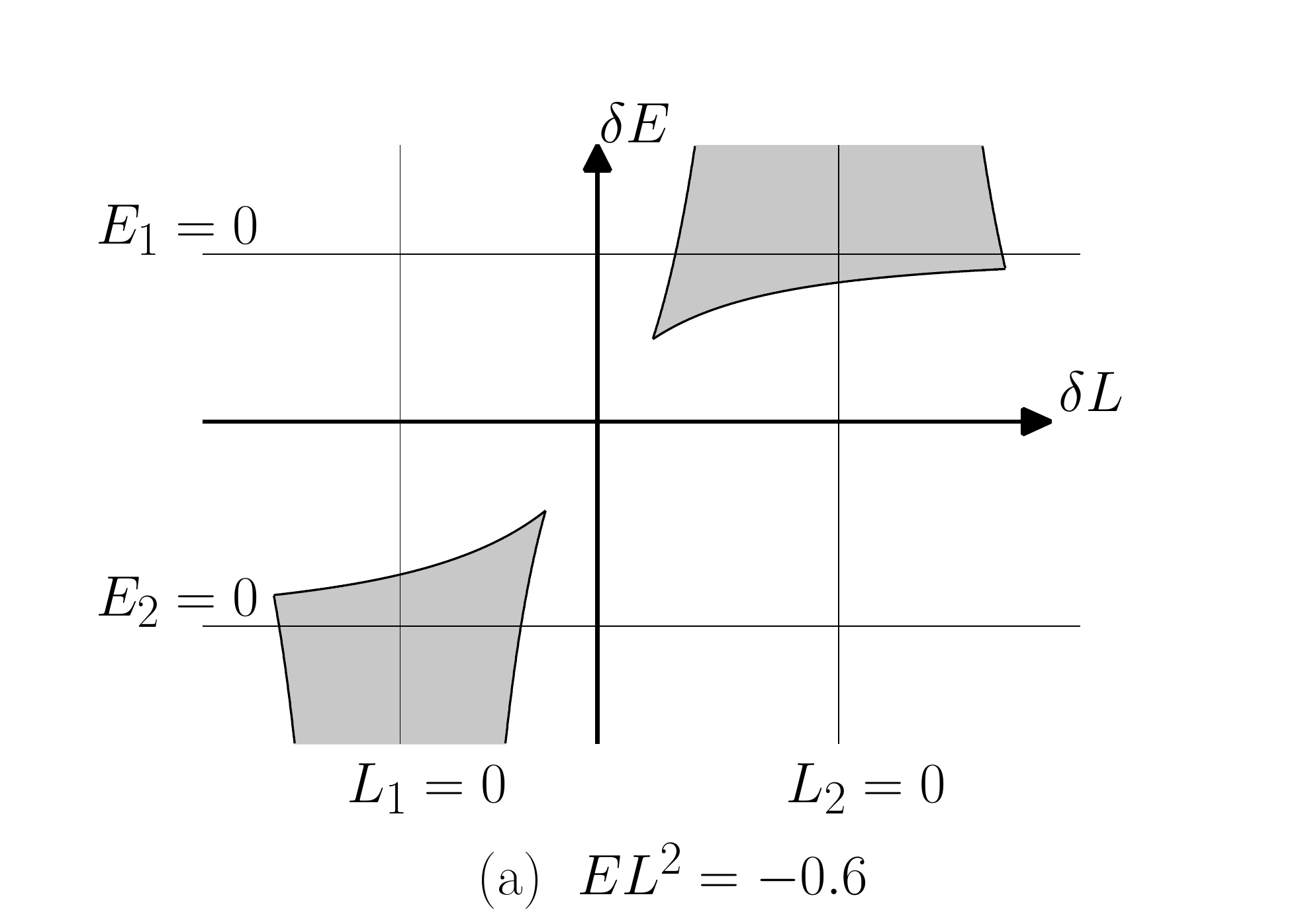} 
  \includegraphics[scale=0.3]{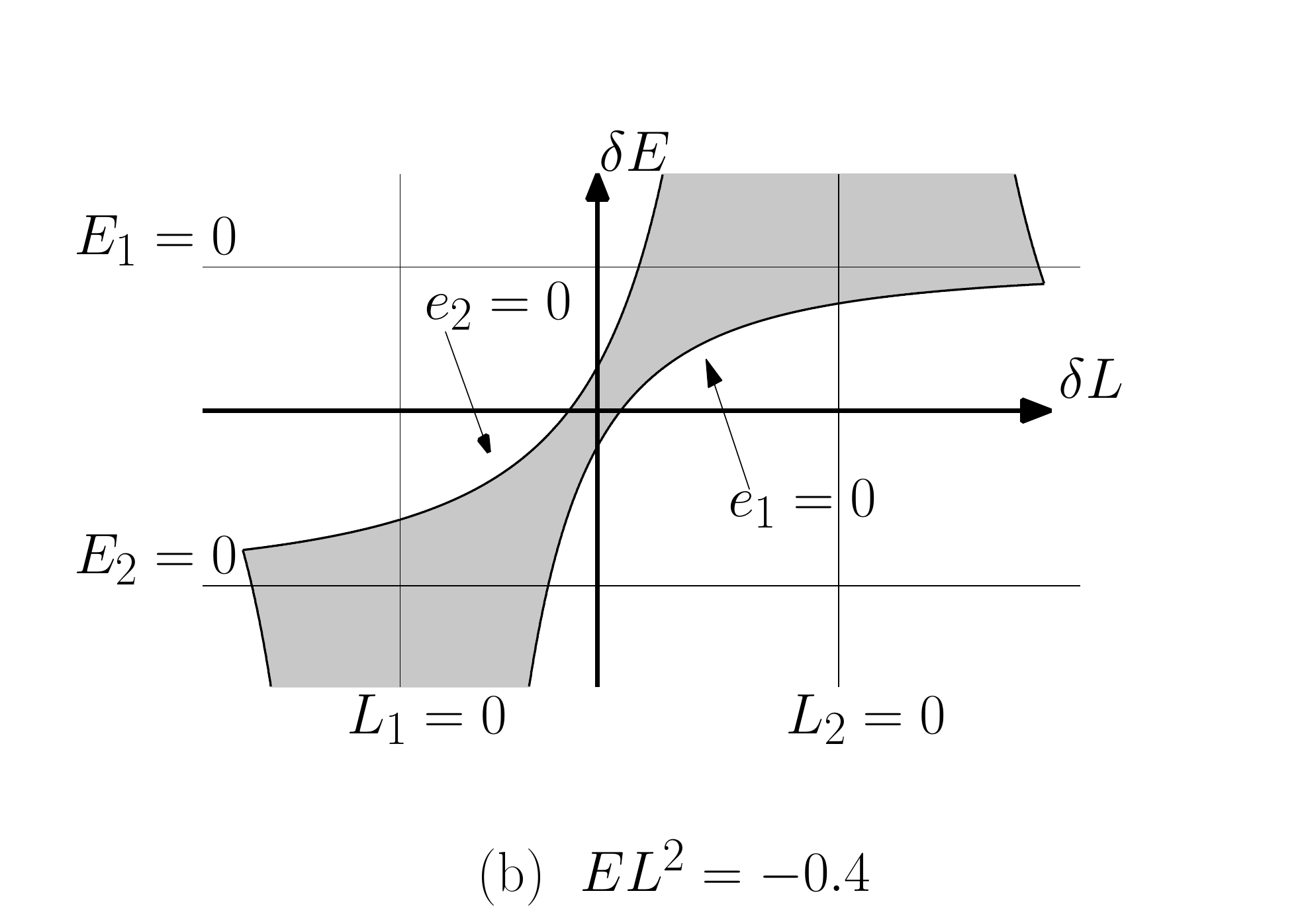} 
  \caption{$\mu_1 = 0.45$: the region of the admissible 
    values for $EL^2 = -0.6$ and $EL^2 = -0.4$.}
  \label{fig:es}
\end{figure*}

Not all the values of $\delta E$ and $\delta L$ correspond
to couple of orbits: namely 
the energy $E_i$ and the angular momentum $L_i$ must satisfy the 
condition $0 \le e_i^2 = 1 + 2  E_iL_i^2$, i.e.
\begin{equation}
\label{esistenza}
\frac 1{\mu_2} \left( -E - 
\frac 1{2L_1^2}\right) \le \delta E \le
\frac 1{\mu_1} \left( E + 
\frac 1{2L_2^2}\right) 
\end{equation}
These inequalities define, in the space $\delta L,\delta E$, the region 
of admissibility
$$A = \{ (\delta L,\delta E)|\, \text{eq.s \eqref{esistenza} hold}\},
$$
which we show in fig. \ref{fig:es}. The boundary of $A$ 
corresponds to $e_1 = 0$, i.e. $\delta E = -\frac 1{\mu_2} \left( E + 
\frac 1{2L_1^2}\right)$, and 
$e_2=0$ i.e. 
$\delta E \frac 1{\mu_1} \left( E + 
\frac 1{2L_2^2}\right)$.
As follows from easy calculations, the topology of the set $A$ 
depends on the value of $EL^2$.
If $EL^2<-1/2$, as in fig. \ref{fig:es} (a), the region is not connected
(this condition is equivalent to the non existence 
of the 'mean orbit', i.e. the orbit of energy $E$ and angular momentum $L$,
whose eccentricity is $\sqrt{1+2EL^2}$). 
If $-1/2\le EL^2 < 0$,  as in fig. \ref{fig:es} (n),  
the region is connected.

If $\delta E \in \left(E/\mu_1,-E/\mu_2\right)$  
the orbits are both elliptic,
while if $\delta E>-E/\mu_2$ (i.e. $E_1>0$) the first orbit is hyperbolic
and if $\delta E<E/\mu_1$ (i.e. $E_2 > 0)$ the second orbit is hyperbolic.
If $\delta L \in \left( -L/\mu_1, L/\mu_2\right)$ both the particles 
move counterclockwise, while if $\delta L = L/\mu_2$
(i.e. $L_1=0$) or $\delta L = -L/\mu_1$ (i.e. $L_2 = 0$)
one of the orbits degenerates.

\vskip3pt
We fix $\mu_1 \le \mu_2$ without loss of generality, 
and we consider the case $E<0$. 
\begin{theorem}
\label{teo:inter}
Let be 
\begin{equation}
\label{sigma}
\sigma= \sigma(\mu_1,\mu_2) = 
- \frac{(1-e^2)(\mu_1^2  + \mu_2^2e)^2}{2\mu_2 e^2}
\end{equation}
where
\begin{equation}
\label{e-sol}
e=\left(\sqrt{(\mu_1/\mu_2)^4 + 8(\mu_1/\mu_2)^2}- 
(\mu_1/\mu_2)^2\right)/4.
\end{equation}
If initially
\begin{equation}
\label{condizione}
EL^2 <  \sigma\\
\end{equation}
the orbits remain elliptic for all times.

Moreover, $|L_1|$ and $|L_2|$ are bounded and 
$e_i \le c_i < 1$ 
for suitable constants $c_1$, $c_2$.
\end{theorem}
\noindent
{\bf Remarks.} 
\begin{enumerate}[{ i.}]
\item As we will see in the proof, if the orbits intersect, 
the value of $EL^2$ is bounded from below:
 $EL^2 \ge - 1/2$.
\item According to the spatial scale invariance of the 
problem, 
the behavior of the system depends only on
the product $EL^2$ of the two invariant quantities $E$ and $L$.
\item In the proof, we study the condition of intersection 
in terms of the orbital parameters; 
for a similar analysis see \cite{laskar0}.
\item If $M$ is the mass of the central body, $G$ the Newton constant
and $k=GM$ the standard gravitational parameter of the system,
the major semi-axes of the orbit $i$ is $L_i^2/(k(1-e_i))$
where the 
eccentricity is  $e_i=\sqrt{ 1 + 2 L_i^2 E_i /k^2}$.
The condition \eqref{condizione} must be rewritten in term of $EL^2/k^2$.
\item The invariant region is large, from two point of view:
it contains couples of orbits which can be 
very different, and the orbits live in a
huge region on the configuration space. For instance, 
in the case $\mu_1 = \mu_2=0.5$, the critical value is $EL^2 = -27/64$,
and, for this value, 
if one of the particles has the same orbit of the Earth, the 
other particle can intersect the Jupiter orbit,
and can arrive at $6.95$ U.A. from the Sun.
If we consider two orbits with $L_1=L_2=L$ and $E_1=E_2=E$, their
eccentricity is $\sqrt{5/32} \approx 0.40$ and the 
ratio between the major and the minor semi-axis is approximately $2.33$. 
\item We have fixed the attractive center, therefore 
we are not considering here a three body problem,
in which we can fix only 
the center of mass.
The planetary three body system 
is usually described as a perturbation of the system obtained
in the canonical heliocentric variables neglecting 
the terms of order $m_1m_2$ (see e.g. \cite{laskar}).
This unperturbed system
is a system of two particles moving independently on 
Keplerian orbits, with respect the position of the large body.
The standard gravitational parameters are $G(M+m_i)$, and can be 
different for the two particles. Our results also hold in
this case, with minor modifications.
\item We have done some preliminary numerical simulations for this model,
with $n>2$ particles.
For $n=3$, if initially all the particles can
collide with the others, in the elastic case
one of the  particle leaves the system after few collisions; 
in the inelastic case  one of the  particle 
stops to interact with the others  after few collisions, 
and the orbits of the others two particles converge.
A system  of a large amount of colliding particles 
exhibits a complex behavior:
a certain number of particles (decreasing with the 
parameter of inelasticity  $\eps$) 
leaves the system, the others particles
asymptotically separate in non interacting 
clusters of one particle or two colliding particles.
Let us note that for some different inelastic model
(somewhat artificial), it can be proved the existence of 
``ringlets'' i.e. the existence of a state of $n$ particle
which does not cease to interact, and whose 
orbits converge (see \cite{kawai}). 

\end{enumerate}

\vskip3pt
\proof

\par
\noindent
We prove the theorem in the case of elastic collisions, for which 
$E$, $L$  and $EL^2$  are conserved quantities, and then the condition 
\eqref{condizione} is invariant for the dynamics.
In the case of inelastic collisions, the thesis  
follows from the fact that $L$ is conserved and $E$ and 
$EL^2$ can only decrease, then the condition \eqref{condizione} 
is invariant for the dynamics also in this case. 

The proof follows from the fact that
if $EL^2$ is sufficiently close
to $-0.5$, and the orbits intersect, 
then the two orbits are elliptic, as we now show.


The two orbits intersect if 
 $\ve x_1 = \ve x_2$ for some value of the anomalies
$\vartheta_1$ and $\vartheta_2$.
Using eq. \eqref{ri}, this condition is expressed by the
equalities
$$\vartheta_2 - \vartheta_1 = \omega_1 - \omega_2 = \delta \omega
\ \ \text{ and } \ 
L_1^2(1 + e_2 \cos \vartheta_2) = L_2^2 (1 + e_1 \cos \vartheta_1).$$
Inserting in the last equation that 
$\vartheta_2 = \vartheta_1 + \delta \omega$,
we obtain an equation in the unknown  $\vartheta_2$ that can be solved if and 
only if
\begin{equation}
\label{interserzione}
e_1 e_2 \cos \delta \omega \le 1 + L_2^2 E_1 + L_1^2 E_2
\end{equation}
Let us define the set of the values of $(\delta L,\delta E)$
for which two orbits of parameters $L_1,E_1$ and $L_2,E_2$ 
intersect if the angle between the periapsides is $\delta \omega = \eta$:
\begin{equation}
\label{def:I}
I_{\eta } = \{ (\delta L,\delta E)\in A |\, 
e_1 e_2 \cos\eta \le 
1 + L_2^2 E_1 + L_1^2 E_2 \}
\end{equation}
By definition
$$I_{\eta_1} \subset I_{\eta_2}\ \ \text{ if } \eta_1 < \eta_2,$$
and in particular $I_{\eta} \subset I_{\pi}$ if $\eta \in [0,\pi]$.
This implies that, if we 
rotate two intersecting orbits, 
in such a way that the two periapsides become in opposition 
($\delta \omega = \pi$), we obtain two orbits which intersect.
The intersection condition is invariant for the dynamics, then
all the values $(\delta L,\delta E)$
during the evolution are in the set
\begin{equation}
\label{Ipi}
I_{\pi}= \{ (\delta L,\delta E)\in A |\, 
e_1 e_2 \ge 
- (1 + L_2^2 E_1 + L_1^2 E_2) \}
\end{equation}
Therefore the set $I_{\pi}$ is invariant
for the dynamics.

Now we show that $I_{\pi}$ is contained
in the region in which both the orbits are elliptic,
if $EL^2$ is sufficiently small.
The set $I_{\pi}$ as defined in \eqref{Ipi} is the union of
the set of  values of $\delta L$ and $\delta E$ in $A$ which solve
\begin{equation}
\label{condizione-quadra}
e_1^2 e_2^2 = (1+2E_1L_1^2) (1+2E_2L_2^2) \ge
(1+L_2^2 E_1 + L_1^2 E_2)^2
\end{equation}
provided that
\begin{equation}
\label{provided}
1+L_2^2 E_1 + L_1^2 E_2 \le  0,
\end{equation}
and the set of values which solve
\begin{equation}
\label{pimezzi}
1+L_2^2 E_1 + L_1^2 E_2 \ge  0,
\end{equation}
Note that this last equation identifies the region $I_{\pi/2}$, i.e. 
the region of intersecting orbits with perpendicular semi-major axis,
and it is equivalent to 
\begin{equation}
\label{pimezzi-risolta}
\delta E ( \mu_1\mu_2 \delta L + (\mu_1 - \mu_2 ) ) \le
1 + 2 E L^2 - (\mu_1 - \mu_2 ) E \delta L,
\end{equation}
Eq. \eqref{condizione-quadra} is equivalent to 
\begin{equation}
\label{condizione-diseq}
\delta E^2 (\mu_1 L_1^2 + \mu_2 L_2^2)^2  - 2 \delta E (L_1^2 - L_2^2) 
( 1 + E(\mu_1L_1^2+\mu_2 L_2^2))  + E^2 (L_1^2-L_2^2)^2 
\le 0
\end{equation}
which can be solved if
\begin{equation}
\label{condizione-discriminante}
(L_1^2-L_2^2)^2 ( 1 + 2E (\mu_1L_1^2+\mu_2 L_2^2)) \ge 0
\end{equation}
If $|L_1|\neq |L_2|$
this condition is equivalent to 
$1+2E(L^2+\mu_1\mu_2 \delta L^2) \ge 0$
then the set $I_{\pi}$ is non void if and only if 
$1+2EL^2 \ge  0$, and it is bounded
by the condition
\begin{equation}
\label{limiti}
\delta L^2 \le \frac { 1 - 2|E| L^2}{2|E|\mu_1\mu_2}
\end{equation}
The boundary of the region identified by the 
inequality \eqref{condizione-diseq} is given by the functions 
\begin{equation}
\label{soluzione}
\delta E = \frac {L_1^2 - L_2^2}{(\mu_1L_1^2 + \mu_2 L_2^2)^2} 
\left( 1 + E(\mu_1L_1^2 + \mu_2 L_2^2) 
\pm \sqrt{ 1 + 2E(\mu_1L_1^2+\mu_2 L_2^2)}\right)
\end{equation}
which have the sign of $L_1^2 - L_2^2$.

\begin{figure*}
  \includegraphics[scale=0.3]{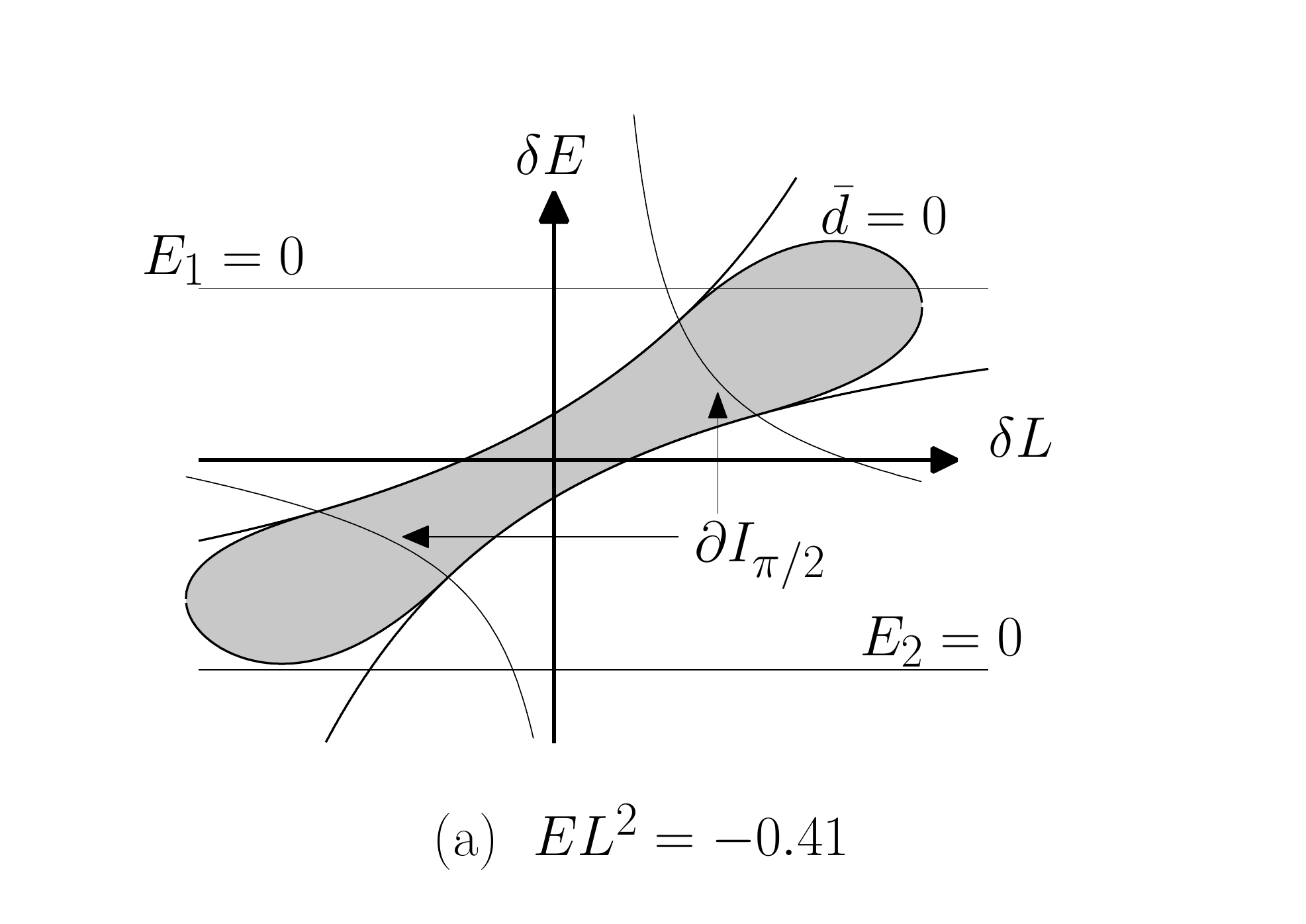} 
  \includegraphics[scale=0.3]{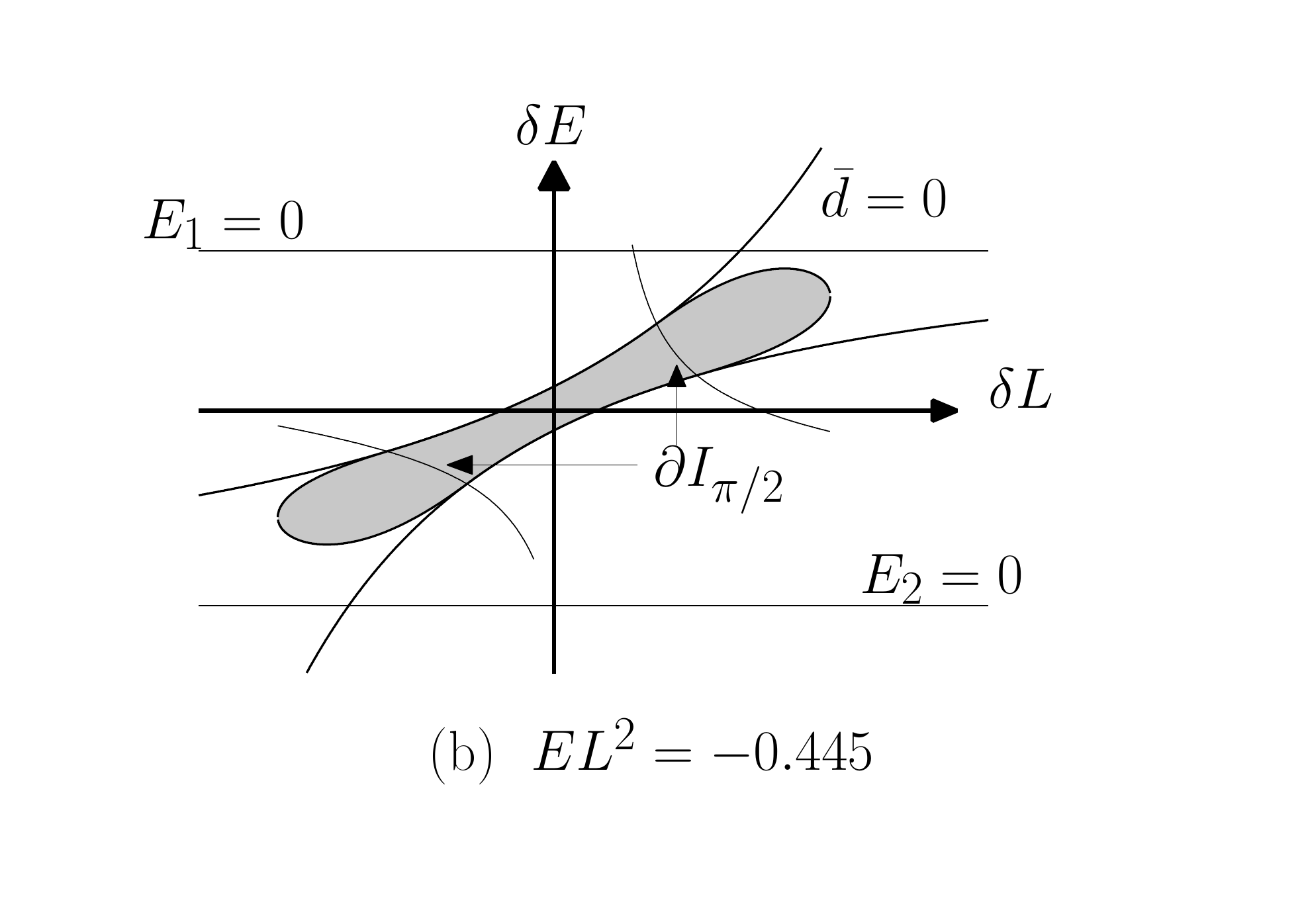} 
  \caption{$\mu_1 = 0.45$: the set $I_{\pi}$ for 
    $EL^2=-0.41$ and $EL^2 = -0.445$.
    Over the line $E_1 =0$,
    the orbit $1$ is hyperbolic. }
  \label{fig:inter}
\end{figure*}

\vskip3pt
In figure \ref{fig:inter} we 
show the region  $I_{\pi}$,  
identified by 
the inequality \eqref{pimezzi-risolta} 
(region $I_{\pi/2} \subset I_{\pi}$) and the inequality
\eqref{condizione-diseq}
(region $I_{\pi}\backslash I_{\pi/2}$).
In figure \ref{fig:inter} (a), the value of 
$EL^2$ is $-0.41$, and the region $I_{\pi}$ 
intersects the  region $E_1 \ge 0$ (i.e. 
$\delta E> |E|/\mu_2$). Then, after a collision, one of the outgoing orbits
can become hyperbolic. In figure \ref{fig:inter} (b),
the value of $EL^2$ is smaller and the invariant
region $I_{\pi}$ is completely contained in the region $\delta E \in 
\left(-|E|/\mu_2,  |E|/\mu_1\right)$ in which both the 
orbits are elliptic.
Then, whatever are the details of 
the collisions, the two orbits remain elliptic,
whit $e_i\le c_i < 1$ for some constants $c_1,c_2$, 
and $|L_1|,\ |L_2|$ are bounded via \eqref{limiti} and \eqref{medie}.

\vskip3pt

Now we will show that the behavior of the system is 
driven by $EL^2$: there exists a critical value which separates 
the two cases.
Let us define 
\begin{equation}
\label{dbar}
\bar d(\delta L, \delta E) = \frac {L_1^2}{1+e_1} - 
\frac {L_2^2}{1-e_2}.
\end{equation}
This quantity is the distance from the periapsis of the orbits
1 and the apoapsis of the orbit 2, in the case of $\omega = \pi$.
The value of $\bar d$ is $0$ on the boundary of $I_{\pi}$
in the first quadrant, out of $I_{\pi/2}$ (see fig. \ref{fig:inter}).
Then the critical value of $E$ e $L$  is such that
$$E_1 = 0,\ e_1 = 1,\ \bar d(\delta E, \delta L) = 0,\ 
\frac {\pa\bar d}{\pa \delta L}  (\delta E, \delta L) = 0$$
(the gradient of $\bar d(\delta E, \delta L)$ is vertical
in the point of tangency to the line $E_1 = 0$).
By deriving $\bar d$ with respect to $\delta L$ we obtain
\begin{equation}
\label{derivata}
\pa_{\delta L} \bar d =
2 \mu_2 \frac {L_1}{1+e_1} +  2 \mu_1 \frac {L_2}{1-e_2} 
- 2 \mu_2 E_1 L_1 \frac {L_1^2}{e_1(1+e_1)^2} +
  2 \mu_1 E_2 L_2 \frac {L_2^2}{e_2(1-e_2)^2} 
\end{equation}
which, using  $-2E_iL_i^2 = 1-e^2_i$, becomes
\begin{equation}
\pa_{\delta L} \bar d = 
2 \mu_2 \frac {L_1}{1+e_1} + 2 \mu_1 \frac {L_2}{1-e_2} 
\mu_2 L_1 \frac {1-e_1^2}{e_1(1+e_1)^2} 
- \mu_1 L_2 \frac {1-e_2^2}{e_2(1-e_2)^2} 
= \frac {\mu_2}{e_1} L_1 - \frac{\mu_1}{e_2} L_2.
\end{equation}
The condition  $\pa_{\delta L} \bar d=0$ and the definition of $L$ 
in  eq. \eqref{medie} allow us to calculate $L_1$, $L_2$
in terms of $e_1$, $e_2$:
\begin{equation}
\label{valoriL}
\begin{array}{l}
L_1 = \mu_1 e_1 L/(\mu_1^2 e_1 + \mu_2^2 e_2)\\ 
L_2 = \mu_2 e_2 L/(\mu_1^2 e_1 + \mu_2^2 e_2)
\end{array}
\end{equation}
Using these expressions in $\bar d = 0$ with $e_1=1$ 
we obtain the following equation
for $e_2$:
$$2\mu_2^2e_2^2 = \mu_1^2 (1-e_2),$$ 
which has only one solution in $(0,1)$, given 
by eq. \eqref{e-sol}.
Substituting this value in the expression of $L_2$,
we obtain the critical values of $EL^2$ as in \eqref{condizione},  
imposing $1+ E_2  L_2^2 =  e_2^2$, with $E_2 = E/\mu_1$ and $e_1=1$:
$$EL^2 = - (1-e_2^2)(\mu_1^2 + \mu_2^2e_2)^2/(2\mu_2 e_2^2).$$
Note that we can find a similar condition for which
the region $I_{\pi}$ is tangent to the line $E_2 = 0$,
but in the case $\mu_1 \le \mu_2$ this second critical value 
of $EL^2$ is greater then the previous, then it can be ignored.
\qed

\section{The invariant region for 
two bodies in $\R^2$}
\label{sezione:estensioni}

In this section, the 
case of two bodies and the case 
of point particles which interact by means of a short range potential are analyzed.

\vskip3pt
\begin{theorem}
\label{teo:D}
We consider two circular bodies in $\R^2$ 
of radii  $R_1$ and $R_2$, 
interacting by means of elastic or inelastic collisions.
If $EL^2 < \sigma(\mu_1,\mu_2)$, with $\sigma$ as defined in 
eq. \eqref{sigma},
and $D=R_1+R_2$ is sufficiently small, then the two bodies remain on 
elliptic orbits.
\end{theorem}

\proof

\par\noindent
We show that if two orbits have points whose distance
is less then or equal to $D$, and 
$EL^2<\sigma$ and  $D$ is sufficiently small, then 
the orbits are elliptic.
We remark that the condition on the distance is preserved 
by the dynamics.

Fixed $L_1,L_2,E_1,E_2$, assuming that the two bodies can collide, 
we need to distinguish two situations:
it exists $\delta \omega$ such that
the orbits intersect, and then $(\delta L,\delta E)\in I_{\pi}$,
or the orbits do not intersect for any $\delta \omega$. 
One of them, named orbit $2$,
is contained in the other.
But
$\min_{\delta \omega} \min_{\vartheta_1,\vartheta_2} |\ve x_1 - \ve x_2|\le D,$
and the minimum is reached for $\delta \omega=\pi$, then 
$0 < \bar d(\delta L,\delta E) \le D$
(note that
if the orbit of the particle $1$ is contained in the orbit of
the particle $2$, 
we have to define $\bar d$ as in \eqref{dbar}, but exchanging the indexes
$1\leftrightarrow 2$).

If $EL^2< \sigma$, the distance between  
the level set $\bar d = 0$ (the boundary of $I_{\pi}$) and 
the critical lines $E_1=0$, $E_2=0$, is strictly positive.
Then, if $D$ is sufficiently small, the set
$\bar d\le D$ does not intersect the region in which 
$E_1\ge 0$ or $E_2\ge 0$,
and this proves the theorem. \qed

\vskip.3cm

\begin{figure*}
  \includegraphics[scale=0.3]{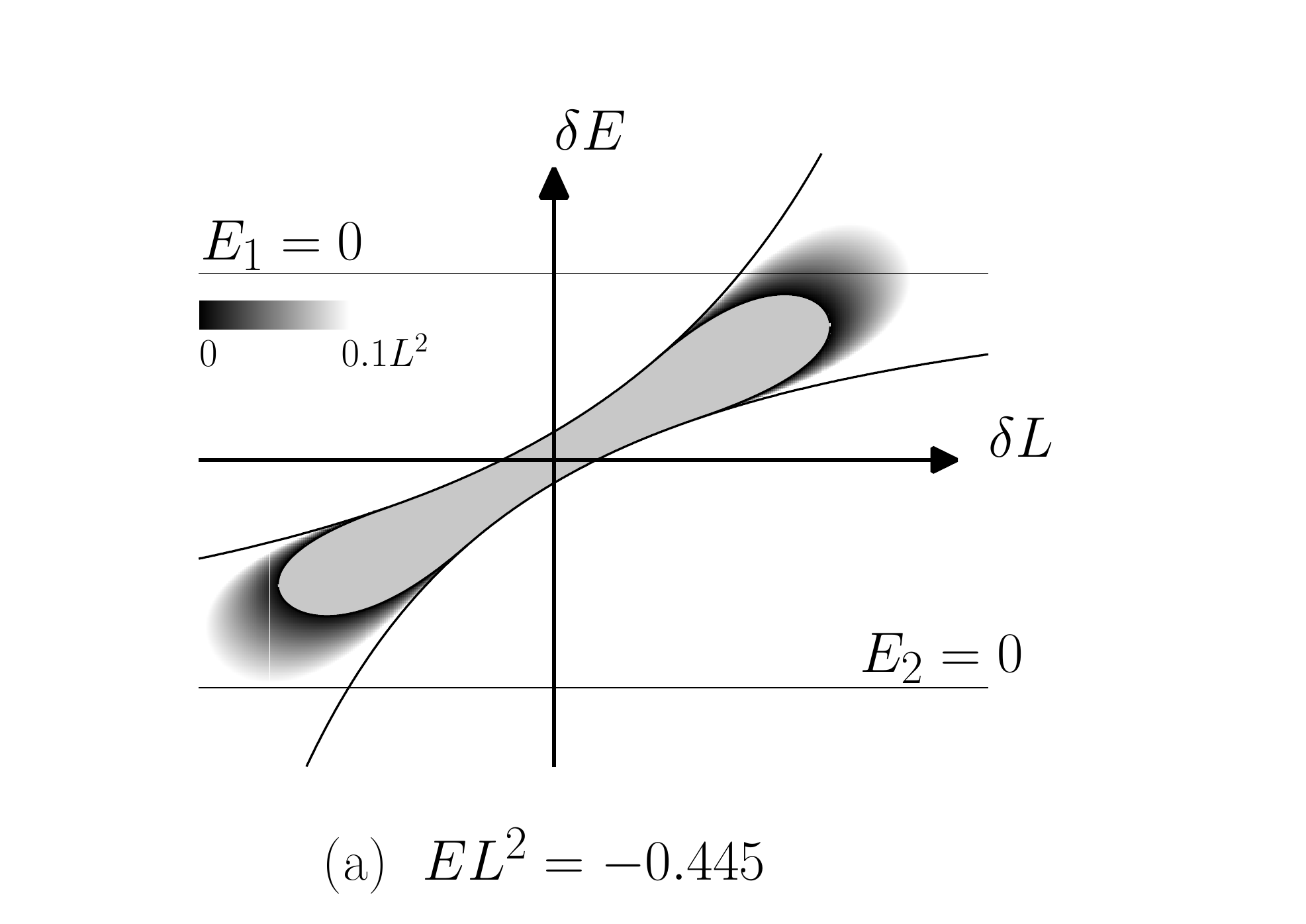} 
  \includegraphics[scale=0.3]{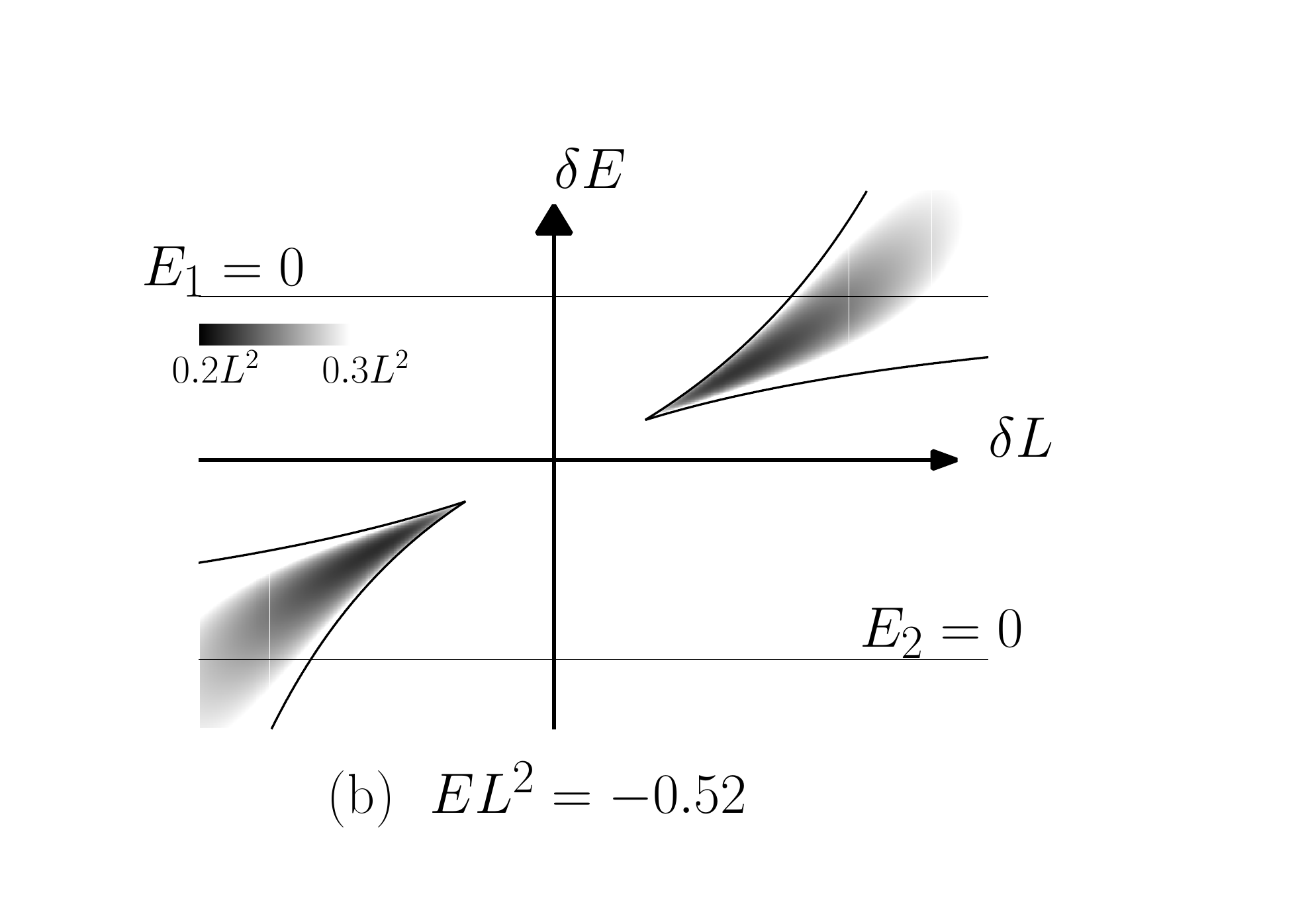} 
  \caption{$\mu_1 = 0.45$: the values of $\bar d$
    in grayscale, 
    in the case $EL^2= -0.445$  
    and $EL^2=-0.52$.
  }
  \label{fig:lev}
\end{figure*}

In figure \ref{fig:lev} we show 
the level sets of $\bar d$.
In figure \ref{fig:lev} (a), $EL^2 = -0.445$ 
and the invariant set $I_{\pi}$ and the values of 
$\bar d$ in the complementary region are shown: the black color corresponds
to $\bar d = 0$, the
white color corresponds to
$\bar d \ge L^2/10$, while the grays correspond to the values  
$\bar d \in (0,L^2/10)$. 
The critical value of $D$ is approximately $0.034 L^2$.
Figure \ref{fig:lev} (b)
is relative to the case $EL^2 = -0.52$, in which the
set $I_{\pi}$ is void. 
Nevertheless, the set $\bar d \le D$
is invariant and is contained in the region of elliptic orbits
if $D$ is sufficiently small. In the graphics,
the black color corresponds
to $\bar d = 0.2 L^2$, the white to $\bar d \ge 0.3 L^2$.
In this case, 
the critical value for $D$ is approximately $0.25 L^2$.

\vskip3pt
We are not able to give a simple expression for the critical values of  
$EL^2$ and $D$, 
but for the case $\mu_1 = \mu_2 = 1/2$, in which 
the particles have the same mass.

\begin{theorem}
\label{teo:D2}
If $\mu_1 = \mu_2 = 1/2$, the conditions on $EL^2$ and $D$ are
$$EL^2 < - \frac{(1-e^2)(1+e)^2}{16e^2}\ \ \text{ where } 
e = \gamma - \sqrt{\gamma^2 - \gamma +1}\ \ \text{ with } \gamma = 2L^2/D > 1
$$
\end{theorem}

\proof

\par\noindent
We proceed as in the proof of Theorem \ref{teo:inter}.
Assuming $\mu_1 = \mu_2 = 1/2$ and $EL^2<-27/64$,
the critical condition $\pa_{\delta L}  \bar d(\delta E, \delta L) = 0$
allows us to obtain the value of $L_1$ and $L_2$
as in eq. \eqref{valoriL}, which becomes 
\begin{equation}
\label{valori2L}
\begin{array}{l}
L_1 = 2e_1 L/(e_1 + e_2)\\ 
L_2 = 2e_2 L/(e_1 + e_2)
\end{array}
\end{equation}
Using these values  in $\bar d = D$ with $e_1=1$,
we obtain the following equation for $e_2$
$$2L^2 (1-2e_2) = D (1-e_2^2),$$
which is solved in $(0,1)$ by 
$$e_2 = \gamma - \sqrt{\gamma^2 - \gamma + 1}
\ \ \text{ where } \ \ \gamma = 2L^2/D
\ \ \text{ with } \ \ D<2L^2.
$$
The corresponding value of $EL^2$
is given by 
$EL^2 = - (1-e_2^2)(1+e_2)^2/(16e_2^2)$.\qed

\vskip.3cm

\begin{figure}
  \includegraphics[scale=0.7]{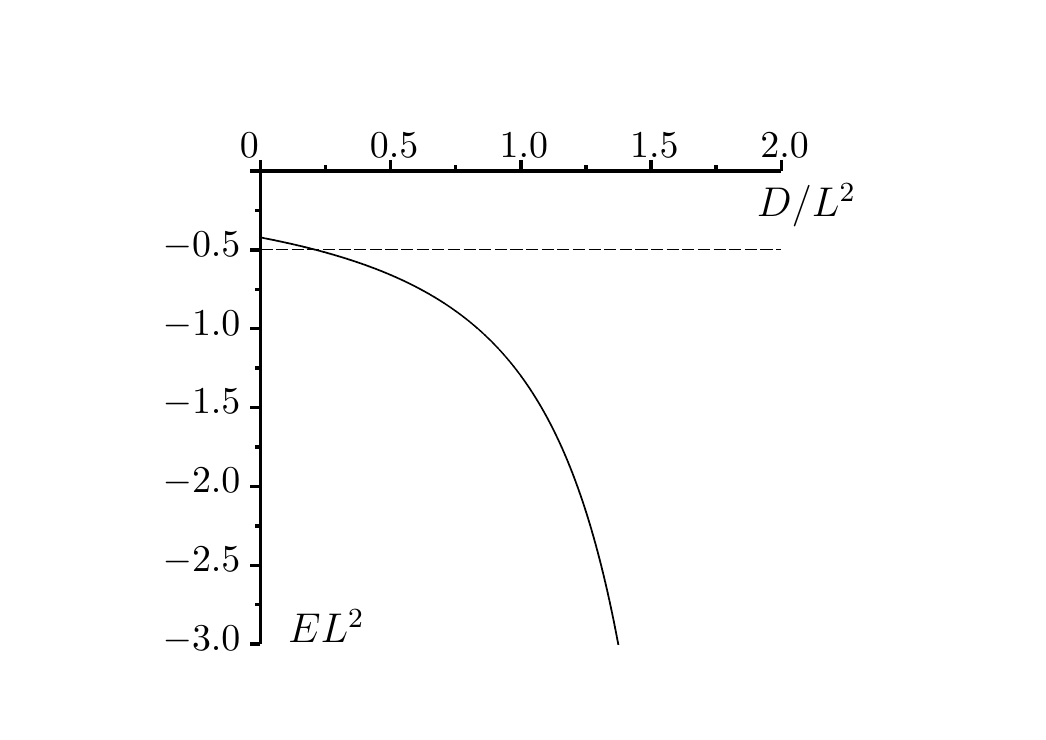} 
  \caption{$\mu_1 = \mu_2 = 0.5$: the critical value of
    $EL^2$ in function of $D/L^2$.
  }
  \label{fig:funz}
\end{figure}

In figure \ref{fig:funz} we plot these values in function
of $D/L^2$. Let us recall that for $EL^2<-0.5$ the 
two orbits do not intersect.

\vskip3pt
The last extension we consider in the two dimensional case is that of two 
point particles interacting by means of a force of symmetric potential 
energy $V(|\ve x_1 - \ve x_2|)$, with a compact support.

\begin{theorem}
\label{teo:V}
Let us consider two point particles interacting by means of a force
of potential energy $V=V(|\ve x_1 - \ve x_2|)$, such that 
$V(r) = 0$ if $r\ge D$ for some $D>0$
and $V(r)\ge -U$ with $U\ge 0$.

If $EL^2 < \sigma(\mu_1,\mu_2)$, with $\sigma$ as 
in eq. \eqref{sigma}, and
and $D$ and $U$ are sufficiently small, 
then the particles remains on bounded orbits.
\end{theorem}

\proof

The specific 
angular momentum $L=\mu_1 L_1 + \mu_2 L_2$ is conserved 
at any time also in this case, because the potential energy is symmetric. 
If $|\ve x_1 - \ve x_2| \ge D$ the energy of the interaction is
zero and the specific energy of the system is exactly 
$$E= \mu_1 \left( \frac {\ve v_1^2}2 - \frac1{|\ve x_1|} \right) + 
\mu_2 \left( \frac {\ve v_2^2}2 - \frac1{|\ve x_2|} \right)$$
Therefore, 
when the two particles leave the region of
the interaction, we can apply Theorem \ref{teo:D}, 
concluding  that the particles remain on elliptic orbits until the
next interaction.

To achieve the proof, we have to discuss 
the motion of the particles during the interaction, 
i.e. when their distance is less than $D$.
If $|\ve x_1 - \ve x_2| < D$,  we can consider
the two 'osculating' Keplerian orbits, i.e. the Keplerian orbits
which correspond to the two couple position-velocity $(\ve x_1,\ve v_1)$,  
and $(\ve x_2,\ve v_2)$.
The specific energies of this two orbits 
are
$$E_1 = \ve v_1^2/2 - 1/|\ve x_1|,\ \ \ \ E_2 = \ve v_2^2/2 - 1/|\ve x_2|.$$
These quantities are not the specific energies of the two particles
(because the contribution of the interaction is not zero), but verify
$$\tilde E = \mu_1 E_1 + \mu_2 E_2 = E - V(|\ve x_1 - \ve x_2|)/(m_1+m_2)
\le E + U/(m_1+m_2)
$$
by the conservation of the total energy.
If $EL^2 < \sigma$ and $U$ is sufficiently small, we have also that
$\tilde E L^2 < \sigma$, then 
we can apply Theorem \ref{teo:D}, using $\tilde E$ instead of $E$.
Therefore, if $D$ is sufficiently small,
as in the hypothesis of Theorem \ref{teo:D}, 
the two osculating orbits are elliptic, with bounded value of $L_1$, $L_2$,
and $e_i\le c_i < 1$, as follows from the compactness
of the set 
$I_{\pi} \cup \{ \bar d(\delta L,\delta E)  \le D\}
\subset \left(-L/\mu_1,L/\mu_2 \right) 
\times \left(E/\mu_1, -E/\mu_2\right)$.\qed

\section{The invariant region for two bodies in $\R^3$}
\label{sezione:dim3}
Here we discuss
the three dimensional case.
We indicate with $\ve L_i=\ve x_i \wedge \ve v_i$ 
the vector which express the specific angular momentum of the 
particle $i$ in the positi`on $\ve x_i$ whit velocity $\ve v_i$, 
and with $\ve L$ 
the specific angular momentum of the whole system
$\ve L = \mu_1 \ve L_1 + \mu_2 \ve L_2$, which is a conserved vector.
In this case, they hold the analogous of 
Theorems \ref{teo:inter}, \ref{teo:D}, \ref{teo:V}, 
where the role of $L$ is played by $|\ve L|$. We summarize these results
in the following theorem.

\begin{theorem}
\label{teo:3d} \phantom{nullaaaaa}

\begin{enumerate}
\item The invariant region for the orbital dynamics in $\R^3$, defined 
as in section \ref{sezione:invariante}, is given by 
$$E|\ve L|^2 <  \sigma$$
with $\sigma=\sigma(\mu_1,\mu_2)$ defined as in eq. \eqref{sigma}.
\item The invariant region for the collisional dynamics of  
two hard spheres in $\R^3$, of radii $R_1$ and $R_2$, with 
$D=R_1+R_2$,
is given by $E|\ve L|^2 < \sigma$
with $D$ sufficiently small.
\item The invariant region for two point particles interacting by means a 
potential energy $V$ 
as in Theorem \ref{teo:V}, is given by  $E|\ve L|^2<\sigma$ with $D$ and 
$U$ sufficiently small.
\end{enumerate}
\end{theorem}
\proof

We prove the theorem staring from the last case, which includes the others
as particular ones.
As in the proof of Theorem \ref{teo:V}, we define
$$E_1= 1/2 \ve v_1^2 - 1/|\ve x_1|, \ \ \ 
E_2= 1/2 \ve v_2^2 - 1/|\ve x_2|, \ \ \ 
\tilde E = \mu_1 E_1 + \mu_2 E_2$$
and we note that
$$\begin{array}{lll}
\tilde E = E &\ \ \ \  \text{ if }\ \ \ \  & |\ve x_1 - \ve x_2| \ge  D \\
\tilde E \le E + U/(m_1+m_2)  &\ \ \ \   \text{ if } \ \ \ \ 
& |\ve x_1 - \ve x_2| <  D.
\end{array}
$$
We indicate whit $\tilde L_i = |\ve L_i|$ the modulus of the 
angular momentum of the orbit of the particle $i$, and we define
$$\tilde L = \mu_1 \tilde L_1 + \mu_2 \tilde L_2$$
which verifies 
$$|\ve L| = |\mu_1 \ve L_1 + \mu_2 \ve L_2| \le \tilde L
\ \ 
\text{ and } \ \ 
E\tilde L^2 \le E|\ve L|^2$$
($E$ is negative).
Moreover
$$\tilde E \tilde L^2 \le \left(E  + \frac{U}{m_1+m_2} \right) |\tilde L^2|$$
then, if $E|L^2| < \sigma$ as in the hypothesis, 
for $U$ sufficiently small,
it also holds
\begin{equation}
\label{condtilde}
\tilde E \tilde L^2 < \sigma
\end{equation}
Let us indicate whit $o_i$ the Keplerian orbit
identified by the position $\ve x_i$ and the velocity $\ve v_i$;
its energy is $E_i$ and its angular momentum is $\ve L_i$.
If the particles can interact, $o_1$ and $o_2$ have
points at distance less then $D$.
We consider the orbits $\tilde o_1$ and $\tilde o_2$ we
obtain rigid rotating in $\R^3$, around the center, 
$o_1$ and $o_2$, in such that 
$\tilde o_1$ and $\tilde o_2$ are in the same plane,
and  the periapsides of $\tilde o_1$ and $\tilde o_2$
are in opposition.
The energy and the eccentricity of the orbits $\tilde o_i$ are
the same of the orbits $o_i$, while we can identify the 
angular momentum of $\tilde o_i$ with the positive scalar 
quantity $\tilde L_i$. 
This two planar orbits intersect or have points at distance less that
$D$. In the first case, we can apply Theorem \ref{teo:inter} 
using \eqref{condtilde}, and we can conclude that the $\tilde o_1$,
$\tilde o_2$, and then $o_1$ and $o_2$ are elliptic.  
In the second case, we can apply Theorem \ref{teo:D}, and we can again 
conclude that, if $D$ is sufficiently small, the orbits are elliptic.

The case of spherical bodies can be considered as a particular
case, in which $V=+\infty$, if $|\ve x_1 - \ve x_2|\le D$. Now $U=0$
and we have only to require a sufficiently small
value of $D$. Finally, the case of the orbital dynamics 
can be considered as the particular case in which 
$D=0$. The hypothesis
$E|\ve L|^2<\sigma$ is then sufficient to achieve the thesis. \qed

\vskip3pt
\noindent
{\bf Remarks.}
\begin{enumerate}[{i.}]
\item The values of $D$ in Theorem \ref{teo:D} 
can be very large with respect 
to the scale $L^2$ of the semi-axis of the orbits, but, if it is so, 
the system is in the region $EL^2 < -1/2$,
in which the orbits do not intersect. Therefore, the system is made of 
two particle which can interact only by means of 
grazing  collisions.

\item It can be interesting to analyze the 
case of Theorem \ref{teo:V} when $V$ is unbounded from below.
The particles can leave the 
system only if their distance remains less then $D$, and, in this case,
we can expect that the center of mass of the two particles
moves on an approximately elliptic orbit. 
On the other hand, there are no a priori bounds on the 
kinetic energy or on the position of this center of mass.

\item It can be also interesting to analyze the case of $n>2$  
particles, with positive radii. In particular it can be expected that
there exist stable ringlets in which the collisions are grazing.
\end{enumerate}

\vskip.5cm

\noindent{\bf Acknowledgments}

\thanks{The authors thank  L. Biasco,  E. Caglioti, G. F. Gronchi, 
P. Negrini, for usefull suggestions on this subject.}

\end{document}